\documentclass[12pt]{article}

\usepackage{graphicx}
\usepackage{amssymb}
\usepackage{amsmath}
\usepackage{color}
\usepackage{algorithm2e}
\usepackage{subcaption}
\usepackage{ulem}

\newtheorem{remark}{Remark}

\newcommand{\kom}[1]{}
\renewcommand{\kom}[1]{{\bf [#1]}}
\definecolor{darkjunglegreen}{rgb}{0.1, 0.14, 0.13}
\definecolor{blau}{rgb}{0.1,0.0,0.9}

\newcounter{komcounter}
\numberwithin{komcounter}{section}

\begin{document}

\title{\bf Management strategies for run-of-river hydropower plants \\
- an optimal switching approach}
\author{%
Niklas L.P. Lundstr\"{o}m\thanks{Department of Mathematics and Mathematical Statistics, Ume\aa\ University,
SE-90187 Ume\aa , Sweden. E-mail: niklas.lundstrom@umu.se} \and
Marcus Olofsson\thanks{Department of Mathematics, Uppsala University, SE-75106 Uppsala, Sweden. E-mail: marcus.olofsson@math.uu.se} \and
Thomas \"{O}nskog\thanks{Statistics Sweden, Stockholm, SE- 171 54 Solna, Sweden. E-mail: thomas.onskog@scb.se}}

\date{}

\maketitle

\begin{abstract}
The mathematical theory for optimal switching is by now relatively well developed, but the number of concrete applications of this theoretical framework remains few. In this paper, we bridge parts of this gap by applying optimal switching theory to a set of production planning problems related to hydropower plants. In particular, we study two different cases involving small run-of-river hydropower plants and show how optimal switching can be used to create fully automatic production schemes in these cases, with non-zero cost of switching between different states of production. Along the way of deriving these schemes, we also create a model for the random flow of water based on stochastic differential equations and fit this model to historical data. This stochastic flow model, which should be of independent interest, mimics the long term seasonal behaviour of the flow while still allowing for stochastic fluctuations and can incorporate a given forecast to damp the impact of such fluctuations in near time. We benchmark the performance of our model using actual flow data from a small river in Sweden and find that our production scheme lies close to the optimal, within 2 \% and 5 \%, respectively, in a long term investigation of the two plants considered.
\newline  \newline
\noindent
\textit{2010 Mathematics Subject Classification: 35K45 
, 49N90 
, 90B50.} 
\\
\noindent
\textit{Keywords:  Stochastic differential equation, river flow model, partial differential equation, variational inequality, production planning}
\end{abstract}


\setcounter{equation}{0} \setcounter{theorem}{0}

\section{Introduction}

Small-scale hydropower plants are in many cases of ``run-of-river" type (ROR) meaning that any dam or barrage is small,
usually just a weir, and generally little or no water can be stored. ROR hydropower plants preserve the natural flow of the river (besides of course at the location of the power plant) and are therefore among the most environmentally benign existing energy technologies. 
Due to a low installation cost, ROR hydropower plants are often cost-efficient and, as such, 
useful both for rural electrification in developing countries and for new hydropower developments in industrialized countries \cite{P02}.
ROR hydropower plants are common in smaller rivers but also exist in larger sizes such as the Niagara Falls hydroelectric plants (Canada/USA), the Chief Joseph dam on the Columbia River (Washington, USA) or the Saint Marys Falls hydropower plant in Sault Sainte Marie (Michigan, USA).

The optimal sizing of ROR hydropower plants has been considered by several authors, see e.g. \cite{AP07, BMM11, F83, GZV09} 
and the references therein. In the estimation of the performance of a given power plant, all these authors omit the cost of switching between different production modes. Doing so, the optimal management strategy can be found trivially by starting the machine when flow is sufficient and stopping it when flow is insufficient. However, in rivers with large and rapid flow fluctuations, which is typically the case in smaller unregulated rivers, such strategies can lead to a large number of switches, the cost of which can not be neglected. For example, starting and stopping the turbines induces wear and tear on the machines and may also require intervention from personnel. Moreover, each start and stop involves a risk which can be considered a cost. To give an example, the major breakdown in the Akkats hydropower plant (Lule river, Sweden) 2002 was caused by a turbine being stopped too quickly, resulting in rushing water destroying the foundation of both the turbine and the generator \cite{AKKATS1, AKKATS2}.

In this paper we create fully automatic production schemes for ROR hydropower plants, with stochastic flow of water and with non-zero cost of switching between different states of production. Our method is based on optimal switching and provides, to the best of our knowledge, a novel way to handle hydropower production planning using stochastic control theory. Along the way of deriving the production strategies, we also create a model for the random flow of water based on stochastic differential equations (SDEs) and fit this model to historical data. The stochastic flow model can incorporate a given flow forecast in order to catch up short term fluctuations.

Our main focus lies on the management of a single power plant without storage capacity or means to control the flow of water, i.e.,~a ROR hydropower plant. Although this is rather restrictive from a practitioner's point of view, we stress that the mathematical framework presented here can be used also for more intricate optimization hydropower production problems, e.g.,~those including dams or pumped-storage hydropower. We have chosen to stay within a simple setup here to highlight the specific features of optimal switching and intend to treat more involved hydropower plants in a series of forthcoming papers. 

The rest of this paper is outlined as follows. In Section \ref{sec:litreview} we give a literature survey and explain the contribution of this paper. Section \ref{sec:optimalsw} contains an introduction to optimal switching problems and, in particular, an outline of the theory in the context of hydropower planning. Sections \ref{sec:flow} and \ref{sec:powerplants} contain models for the flow of water and the power plants under consideration, respectively. Section \ref{sec:numericsPDE} contains a very brief outline of the numerical approach taken to solve the variational inequalities appearing. We thereafter present the result of our parameter estimation and the performance of the constructed strategy in Section \ref{sec:results}. We end with Section \ref{sec:discussion} in which we discuss our approach in general, and our results in particular, and make some concluding remarks.


\setcounter{equation}{0} \setcounter{theorem}{0}

\subsection{Literature survey and our contribution} \label{sec:litreview}

Optimal switching is a relatively new and fast growing field of mathematics combining optimization, SDEs and partial differential equations (PDEs) \cite{BGSS20, BJK10,DHP10, AF12, AH09, HM12, HT07, K14, LNO14, LNO14b, LOO17, M16_1, M16_2, P18, P20, P20_2}. However, a literature survey shows that,
although the mathematical theory is well developed, applications of optimal switching to real life problems is a far less explored area. A possible explanation for this discrepancy could be the difficulty of formulating real problems in mathematical terms and, conversely, to interpret the theoretical results in practical terms. 

Most commonly, applications are found in the context of real options, see \cite{BO94, BS85, CL08, DZ01, GKL17}. In \cite{ACLP12} the authors provide a probabilistic numerical scheme for optimal switching problems and test their scheme on a fictitious problem of investment in electricity generation. In \cite{CL10} valuation of energy storage is studied with the help of optimal switching and in \cite{P20_3, PE17} the authors study how the framework can be used to track electricity demand. 


This paper extends the use of optimal switching by applying the general framework described in Section \ref{sec:optimalsw} to two canonical examples of ROR hydropower plants, consisting of one and two adjustable units, respectively, more thoroughly described in Section \ref{sec:powerplants}. 
In particular, we construct a management strategy by solving a variational inequality related to the optimal switching problem (see \eqref{pde}). However, for this solution to be practically meaningful, we first adapt and calibrate the optimal switching problem to the case of ROR power plants.

A main feature of optimal switching-based production planning is that it allows for random factors influencing the production strategy (for details see Section \ref{sec:optimalsw} below). This randomness is in general given by a Markovian stochastic process and can incorporate any number of different variables. However, to reduce computational strain we will in this paper focus on how the flow of water impacts the production. 

Popular streamflow models include linear time series models, such as ARMA models with Gaussian or GARCH noise \cite{MO13, MNK90, WGVM05}, and non-linear time series models, such as SETAR models \cite{FFODN19}. Modern approaches also include neural networks and machine learning techniques, see, e.g.,~\cite{MBFZ17, MS17} and the references therein. Another suitable approach, and the approach that we have adopted here due to its natural relation to the optimal switching framework, is SDEs driven by Gaussian white noise and/or compound Poissonian impulses \cite{BTU87}. 
In particular, we develop a new stand-alone SDE-based model for the flow of water $Q$, based on historical data and driven by Gaussian white noise, which mimics the long term seasonal variations of the flow while still allowing for short term fluctuations and flow forecasts.

When trying to maximize monetary profit, the electricity spot price $P$ at which the electricity produced is sold is of course also of interest when planning the production. In general, $P$ is an exogenous stochastic process which, in principle, can be incorporated into our model similarly to the water flow by constructing an SDE for $P$ and applying optimal switching theory. However, modelling electricity prices is a non-trivial task and prices are usually not It\^o diffusions, but rather discontinuous jump-processes, see, e.g., \cite{HL12, WBT04, WSW04}, making the operator in the variational inequality to be solved non-local (see, e.g., \cite{LNO14, LNO14b}). For simplicity, we therefore let $P$ be a continuous \textit{deterministic} process in the current paper. 
We stress, however, that our approach readily can be extended to random electricity prices (and streamflow models with compound Poisson impulses) at the cost of computational complexity. 

We have chosen not to let deeper knowledge of different types of generators and how they operate be a prerequisite for the understanding of this paper and we hence avoid going into any such technical details. A comprehensive survey of different turbines and generators and their distinct characteristics can be found in, e.g., \cite{B15, B01, W84}.

\section{Optimal switching in the context of hydropower}  \label{sec:optimalsw}

Tailored to the setting of hydropower plants, the optimal switching problem can be described as follows. We consider a manager of a hydropower plant with several units, each unit being a sub-power plant, i.e.,~a turbine and a generator. Each unit can be started and stopped separately in order to adjust the production of the power plant to the supply of water or to the production demand. This implies that the manager has the option to run the plant in $m \geq 2$ production modes, corresponding to running different combinations of units. Starting and stopping units induces wear and tear on the units and therefore the manager finds herself in a trade-off, weighing the benefits of changing production state against the costs induced by making these changes.

Let $X = \{X_t\}_{t\geq 0}$ denote a Markovian stochastic process representing the features which influence the production. For small hydropower plants, $X_t$ may represent the flow of the river, but it may also be interpreted as, e.g.,~production demand for a frequency regulating plant, the spot price of electricity, or a varying cost of production such as for an oil driven power plant. The process $X$ may be multi-dimensional and hence incorporate all of the above and more, but we will in this paper consider $X =(Q,P)$ where $Q$ and $P$ are one-dimensional processes representing flow of water and spot price of electricity, respectively. Moreover, we let $f_i(X_t,t)$ denote the instantaneous payoff generated in production mode $i$ at time $t$, when the state of the underlying process is $X_t$. Depending on the interpretation of $X$ and the choice of $f$, $f_i(X_t,t)$ can be interpreted as, e.g., the power delivered by a power plant or as the instantaneous monetary profit per unit time. Finally, we associate to each start and stop of a unit a cost $c_{ij}$ for switching from production mode $i$ to production mode $j$, where $i,j \in \{1,\dots m\}$.

The manager of the power plant controls the production by choosing a \textit{management strategy}, i.e.,~a combination of a non-decreasing sequence of stopping times  $\{\tau_k\}_{k\geq 0}$, where, at time $\tau_k$, the manager decides to switch the production from its current mode to another, and a sequence of indicators  $\{\xi_k\}_{k\geq 0}$, taking values in $\{1,\dots,m\}$, indicating the mode to which the production is switched. More precisely, at $\tau_k$ the production is switched from mode $\xi_{k-1}$ to $\xi_k$ and when starting in mode $i$ at time $t$, we have $\tau_0 = t$ and $\xi_0 = i$.  We stress that $\tau_i$ is required to be a \textit{stopping time} and as such it is adapted to the filtration $\mathcal F^X$ generated by the underlying process $X$. In less mathematical terms, this simply means that the decision to switch at time $t$ must be based solely on the information made available up to time $t$, i.e., the manager cannot``peek into the future'' when making her decision\footnote{Although this is obviously impossible from a practical point of view, it must be stated as an explicit restriction in the mathematical formulation.}.

A strategy $(\{\tau_k\}_{k\geq 0},\{\xi_k\}_{k\geq 0})$ can be represented by the random function $\mu : [0, T] \to \{1,...,m\}$ defined as
\begin{align*}
\mu_s \equiv \mu(s) = \sum_{k\geq 0} \mathbb{I}_{[\tau_{k}, \tau_{k+1})}(s)\xi_{k},
\end{align*}
which indicates the mode of production at time $s$. Here, $\mathbb{I}_{[\tau_{k}, \tau_{k+1})}(s)$ denotes the indicator function, i.e.,~$\mathbb{I}_{[\tau_{k}, \tau_{k+1})}(s) = 1$ if $s\in [\tau_{k}, \tau_{k+1})$ and 0 otherwise. When the production is run using a strategy $\mu$, defined by $(\{\tau_k\}_{k\geq 0},\{\xi_k\}_{k\geq 0})$, over a finite horizon $[0, T]$, the total expected payoff is
\begin{align*}
E\biggl [\int \limits_0^T  f_{\mu_s}(X_s,s) ds - \sum_{k\geq 1,\tau_k < T} c_{\xi_{k-1},\xi_k} \biggr ].
\end{align*}
Similarly, given that the stochastic process $X$ starts from $x$ at time $t$, the profit made using strategy $\mu$ starting in $\xi_0=i$, over the time horizon $[t, T]$, is
\begin{align} \label{eq:OSP}
J_i(x,t,\mu) = E\biggl [\int \limits_t^T  f_{\mu_s}(X_s,s) ds - \sum_{k\geq 1,\tau_k < T} c_{\xi_{k-1},\xi_k}\bigg| X_t = x \biggr ].
\end{align}

The task of the manager is now to maximize the expected payoff, i.e.~to find the most profitable trade-off between switching to more efficient states and minimizing the total cost of switches. In general, this problem, often referred to as an optimal switching problem, consists in finding the value function
\begin{align*}
u_i(x, t) = \sup_{\mu\in\mathcal A_i} J_i(x, t, \mu),
\end{align*}
where $\mathcal A_i$ is a set of strategies starting from $\xi_0=i$, and the optimal management strategy $\mu^* \in \mathcal A_i$, defined by $(\{\tau^*_k\}_{k\geq 0},\{\xi^*_k\}_{k\geq 0})$,
such that
\begin{align*}
J_i(x, t, \mu^*) \geq J_i(x, t, \mu)
\end{align*}
for any other strategy $\mu \in  \mathcal A_i$. We note that the value function and optimal strategy of an optimal switching problem depends on the set of available modes $\{1,\dots,m\}$, the (finite) time horizon $T$, the running payoff functions $f_i$, $i\in \{1,\dots,m\}$, the switching costs $c_{ij}$, $i,j\in \{1,\dots,m\}$, as well as the dynamics of the underlying stochastic process $X$.

From the dynamic programming principle, see, e.g.,~\cite{TY93}, one can derive that the value function $u(x, t) = (u_1(x, t), \dots, u_m(x, t))$ satisfies the following system of variational inequalities
\begin{align}\label{pde}
&\min \left\{ -\partial _{t}u_{i}\left( x,t\right) - \mathcal{L} u_{i}\left( x,t\right) -f_i\left( x,t\right),
\;u_{i}\left( x,t\right) -\max\limits_{j\neq i}\left\{ u_{j}\left(
x,t\right) -c_{ij}\left( x,t\right) \right\} \right\}  = 0, \\
&u_{i}(x,T) = g_{i}\left( x\right) ,  \notag
\end{align}
for $x\in \mathbb{R}$, $t\in [0,T]$, $i\in \{1,...,m\}$ and where 
$$\mathcal{L} u(x,t) : = \sum_{i=1}^{d_X} b_i(x,t) \frac{\partial }{\partial x_i} u(x,t)+ \sum_{i,j=1}^{d_X} \frac{1}{2} (\sigma \sigma^\ast)_{i,j}(x,t) \frac{\partial^2}{\partial x_i \partial x_j}  u(x,t),
$$ 
is the infinitesimal generator of the underlying $d_X$-dimensional stochastic process $X$ with dynamics
$$
dX_t = b (X_t,t) dt + \sigma(X_t,t) dW_t.
$$
The operator $\mathcal{L}$ is a linear operator whenever the process $X$ is driven by Gaussian white noise (which we assume in this work). 

Moreover, the optimal strategy from state $i$ is given iteratively by the solution to \eqref{pde} by
\begin{align*}
&\tau^\ast_0=0 \qquad \xi^\ast_0=i \\
&\tau^\ast_k = \inf\left \{t > \tau^\ast_{k-1}: u_{\xi^{\ast}_{k-1}}(X_t,t) \leq \max_{j \neq \xi^\ast_{k-1}}\{ u_j(X_t,t) - c_{\xi^\ast_{k-1}j} \} \right \}
\\
& \xi^\ast_k = \underset{j}{\text{argmax}} \left \{ u_{j}(X_{\tau^\ast_{k}},\tau^\ast_{k})  - c_{\xi^\ast_{k-1}j} \right \}.
\end{align*}

Last, we remark that system \eqref{pde} consists of $m$ PDEs with interconnected obstacles and that a unique Lipschitz continuous solution to this system exists under the assumptions of this paper, see, e.g., \cite{ADPS09, AH09, LNO14b}. 

\section{Modeling river flow with an SDE} \label{sec:flow}
In this section we construct an SDE whose solution resembles the actual flow of water in the river under investigation.
This resemblance must of course hold in the long-term (seasonal) sense, but must also allow for short-term fluctuations due to inter-yearly variations. More specifically, we use historical data to find functions $b_s$ and $\sigma_s$, such that the solution to the stochastic differential equation
\begin{align*}
dQ_t =& b_s(Q_t,t)\, dt + \sigma_s(Q_t,t)\, d\tilde{W}_t,\quad \text{if }t\in [s,T],  \notag \\
Q_s=&q
\end{align*}
where $\tilde{W}_t$ is a standard Brownian motion, is similar (in some appropriate sense) to the actual flow of water.

Results indicate that log-transformation of river flow data may increase prediction accuracy of flow models, see \cite{AYKA17}, and we therefore work with the logarithm of the flow rather than the flow directly. We treat the seasonal and short-term resemblance separately, as outlined below.

\subsection{Seasonal variations}

Starting with the long term seasonal variations, we define $R_t = \log Q_t$. We let $r_t$ be a function describing the expected value of the logarithm of the flow at time $t$, independent of current observations. Defined as such, $r_t$ reflects the expected seasonal variation in flow due to spring flood, autumn rains etc., but without any consideration taken to observations from the current year\footnote{As an example, one would expect that heavy snowfall from January to March would increase the probability of a long and intensive spring flood. Such considerations are \textbf{not} included in $r_t$.}, and we may thus estimate the deterministic function $r_t$ from historical flow. More precisely, we will construct $r_t$ as a one-week moving average of the logarithm of the flow, see Section \ref{sec:parametervalues}. The choice of a one-week moving average here is a trade-off between capturing seasonal variations, such as the spring flood, without letting the mean flow depend too much on the flows of particular years.

\subsection{Short term fluctuations}

Next, we consider the fluctuations around the expected yearly mean, $S_t=R_t-r_t$, and assume that these fluctuations are given by an Ornstein-Uhlenbeck process reverting towards $0$, i.e.
\begin{equation*}
dS_t =  -\kappa S_t \, dt + \sigma\, dW_t,
\end{equation*}
where $\kappa>0$ and $\sigma>0$ are constants to be determined and $W_t$ is a standard Brownian motion. By standard It\^o calculus, the flow $Q_t = \exp\left(r_t+S_t\right)$, then satisfies the following stochastic differential equation 
\begin{equation} \label{eq:sde}
dQ_t = \left(r'_t+\frac{1}{2}\sigma^2-\kappa\left(\log Q_t - r_t\right)\right)Q_t\, dt + \sigma Q_t\,dW_t.
\end{equation}
Note that the particular form of \eqref{eq:sde} ensures that $Q_t$ stays positive.

To estimate the parameters $\kappa,\sigma$ we consider the asymptotic variance and asymptotic autocorrelation for lag $\tau$ of an Ornstein-Uhlenbeck process, which are given by
\begin{equation}
\text{Variance}=\frac{\sigma^2}{2\kappa}\qquad\text{and}\qquad\text{Autocorrelation}=e^{-\kappa\tau}.\label{eq:varacf}
\end{equation}
For a given set of historical flow data, we first calculate the logarithm of the data and subtract the running-mean $r_t$ from above to obtain an empirical time series for $S_t$. We calculate the sample autocorrelation function of the time series and estimate the value of $\kappa$ by a linear regression with the logarithm of the sample autocorrelation function as the dependent variable and the lag as the covariate. Finally, we calculate the sample variance of the time series and deduce an estimate of $\sigma$ from the first equation in \eqref{eq:varacf} and the estimated value of $\kappa$.

\subsection{Forecasts}
The SDE constructed so far respects seasonal variations whilst still allowing for inter-yearly forecasts. However, in order to perform optimally, our model must also be able to treat short-term fluctuations based on, e.g.,~weather forecasts or upstream measurements of the flow. Such forecasts will be included in the model by altering the drift in the dynamics of $Q$, i.e.,~by changing \eqref{eq:sde} on a short time span close to the present time $s$. 

More precisely, we let for $t \geq s$
\begin{equation} \label{eq:sde!}
d Q_t = b_s(Q_t,t) \, dt + \sigma_s(Q_t, t) \, dW_t,
\end{equation}
where
\begin{equation}\label{eq:sigma}
\sigma_s(Q, t)   = \sigma Q
\end{equation}
as before and
\begin{align} \label{eq:b}
b_s(Q, t) &=
\begin{cases}
b_s^{f} (t)							 &\textrm{if}\quad s \leq t \leq s+l +\ell,\\
 \left(r'_{t}+\frac{1}{2}\sigma^2-\kappa\left(\log Q- r_{t}\right)\right) Q	& \textrm{if}\quad s+l+\ell < t,
 \end{cases}
\end{align}
where $b^f_s$ is a function constructed from the long term (log-)mean $r$, the flow $Q_s$ at time $s$ and the forecast of the flow during $(s, s+l)$ as follows. In the above, $l$ and $\ell$ are parameters representing the length of the forecast and the (estimated) time it takes for the forecasted flow to return to the long term (log-)mean $r$.
Starting at time $s$ with current flow $Q_s$ and given a forecast $\{F_r\}_{s< r \leq	s+l}$ of the future flow at times $t \in (s,s+l)$ we set
$$
g_s(t) = \begin{cases}
\log(Q_s) & \mbox{for $t=s$} \\
\log (F_t) & \mbox{for $t\in (s,s+l]$} \\
\log(F_{s+l}) +(t-(s+l)) \frac{ r_{s+l+\ell}- \log(F_{s+l}) }{\ell} & \mbox{for $t \in (s+l, s+l+\ell] $}
\end{cases}
$$
and let
$$
b_s^{f} (Q,t)	 =  \left(g_s'({t})+\frac{1}{2}\sigma^2-\kappa\left(\log Q - g_s({t})\right)\right)Q.
$$

More explicitly, we calculate the drift $b^f_s$ as in \eqref{eq:sde} but with the (log-) mean $r$ replaced by $g_s$, where $g_s(t)$ is given directly by the forecast for $t \in (s, s+l]$, coincides with $r_t$ for $t >  s+l+\ell$ and is linearly interpolated between $s+l$ and $s+l+\ell$. The impact of such forecasts are illustrated in Figure \ref{fig:forecasts}.

We stress that as time evolves, the forecast will be updated and the function $b_s^f$ needs to be updated accordingly each time $s$ that a new forecast becomes available. As the starting time $s$ will be clear from context we will drop the subscript $s$ and simply write $b^f(Q,t)$ in the following, although this function varies with $s$ as parameter.


\section{Modeling the payoff structure of power plants} \label{sec:powerplants}
As mentioned above, we will consider two canonical examples of ROR hydropower plants and outline the payoff structure of these in the current section. We will measure the performance in monetary units (m.u.), but one can easily modify the below to have, e.g.,~total electricity produced as the trait for optimization.

\subsection{Power plant I: One adjustable unit}\label{sec:model1}
We consider first the simple case of a hydropower plant having a single unit. The unit is designed for the flow $Q_d$, but can be run over a wide flow range $[Q_{min}, Q_{max}]$ with lower efficiency. We assume that the unit, automatically and at negligible cost, adjusts to the available flow, and the task of the manager is thus to find out when to start and stop the unit. In the setting of optimal switching, we model the above power plant as follows. The power plant can be run in two states, `1' and `0', representing `on' and `off', respectively, and for each switch from $0$ to $1$ or from $1$ to $0$ the manager must pay a cost $c_{01}$ or $c_{10}$, respectively. We assume that the electricity output (in Watts) of the unit when in state $1$ is given by
\begin{equation*}
W_1(Q) =
\begin{cases}
0							& \textrm{if}\quad Q < Q_{min},\\
c\, \eta (Q) \, Q						& \textrm{if}\quad Q_{min} \leq Q < Q_{max},\\
c\, \eta (Q_{max}) Q_{max}					& \textrm{if}\quad Q_{max} \leq Q,
\end{cases}
\end{equation*}
where the constant $c$ is simply given by $c = \rho\, g\, h$, where $\rho = 10^3$ kg/m$^3$ is the (approximate) density of water, $g = 9.82$ m/s$^2$ is the (approximate) gravitational constant, and $h$ is the water head in meters, and $\eta(Q) \in (0,1)$ is the flow-dependent efficiency. In practice, this latter function is given by the specific characteristics of the turbine and generator, but in general it is a concave down function with a maximum at the design flow $Q_d$, see, e.g., \cite[Figure 9]{IPCC}. To mimic such behaviour, we let
\begin{align}\label{eq:efficiency-curve-assumption}
\eta(Q) = \alpha - \beta \left( \frac{Q}{Q_d} - 1\right)^2,
\end{align}
where $\alpha$ is the efficiency at design flow and $\beta > 0$ quantifies the concavity. Multiplying $W_1$ with the spot electricity price and integrating over time gives the income when running the plant in state $1$.

To be able to determine an optimal strategy in terms of monetary profit, we need to consider not only the produced electricity, but also running-costs of the plant as well as electricity prices. 
We assume that the running cost of the unit is $c_{run}$ per unit time, the electricity price at time $t$ is $P_t$, and that a large additional cost $c_{low}$ must be paid for each unit of time that the production unit is run with insufficient flow $Q <Q_{min}$. The cost $c_{low}$ is motivated by excessive deterioration of the power plant and the cost should be high enough to avoid running the unit with insufficient flow. Summarizing the above, we find that a reasonable generic payoff function for a unit is given by
\begin{align}\label{eq:payoff}
f_1(Q,P,t) = - c_{run} +
\begin{cases}
-c_{low}				& \textrm{if}\quad Q < Q_{min},\\
P\,c\, \eta(Q)\,Q 				& \textrm{if}\quad Q_{min} \leq Q < Q_{max},\\
P\,c\, \eta(Q_{max}) Q_{max} 			& \textrm{if}\quad Q_{max} \leq Q,
\end{cases}
\end{align}
where the parameters $c,\eta, c_{low}, c_{run}, Q_{min}$ and $Q_{max}$ in practice are determined by the specific unit under study and where $P$ represents the current spot price of electricity. We assume that all the above numbers are normalized so that we can take $f_0 \equiv 0$, i.e., the running cost/payoff when in state $0$ is $0$.

In the above we have, without loss of generality, omitted fixed costs such as maintenance of buildings, insurances, etc.~as these costs do not influence the optimal management strategy. For simplicity, we have also omitted maintenance of the unit, although planned (or unplanned) production stops for maintenance could potentially influence the optimal strategy. 

\subsection{Power plant II: Two homogenous adjustable units}

In our second case, we expand the above power plant with another identical unit so that there are in total three different states available to the manager; running no unit (state $0$), running a single unit (state $1$), or running both units (state $2$). As the units are assumed to be identical, there is no difference between running unit $1$ or unit $2$ in state $2$\footnote{No difference seen by our model, that is. There could of course be a difference in practice given different wear of the units, planned maintenance, etc.}. We assume that the payoff structure of both units is as specified in \eqref{eq:payoff}.

In this setting, a difficulty beyond that of the extra state is introduced. Indeed, with both units running, the manager may control how much of the current flow $Q_t$ that is directed towards each of the two running units. With $\delta$ representing the fraction of water directed towards unit $1$, we see that the payoff in state $2$ is
$$
f_1(\delta Q,P,t) + f_1((1-\delta)Q,P,t)
$$
and $\delta$ is thus an additional degree of freedom to optimize over. However, as long as the cost of adjusting units and diverting water is negligible, this optimisation can be done independently of the switching problem. The two-dimensional optimization problem involving two homogenous adjustable units can thus be reduced to a pure switching problem by introducing the payoff function
\begin{equation} \label{eq:f3}
f_2(Q,P,t) = \max_{\delta \in [0,1]} \left \{f_1 (\delta Q,P,t) + f_1 ((1-\delta) Q,P,t) \right \}
\end{equation}
in addition to $f_0 \equiv 0$ and $f_1$ as in \eqref{eq:payoff}. The optimization in \eqref{eq:f3} can be carried out analytically if $f_1$ is sufficiently nice. However, since our final scheme is based on numerical solutions to PDEs, it is merely a question of computational power if the optimization problem \eqref{eq:f3} can only be solved numerically. In any case, this optimization can be done in advance without consideration of current and forecasted flow and to arbitrary precision, making the remaining optimization problem purely a question of optimal switching.

\begin{remark}
The above case can easily be expanded to the case of heterogeneous units. In this case, we assume that unit $2$ yields payoff
\begin{align*}
\tilde{f}_1(Q,P,t) = - \tilde c_{run} +
\begin{cases}
-\tilde c_{low}				& \textrm{if}\quad Q < \tilde Q_{min},\\
P\,c\,\tilde \eta(Q)\,Q 				& \textrm{if}\quad \tilde  Q_{min} \leq Q <\tilde Q_{max},\\
P\,c\, \tilde \eta( \tilde Q_{max})\tilde Q_{max} 			& \textrm{if}\quad \tilde Q_{max} \leq Q,
\end{cases}
\end{align*}
i.e.,~unit $2$ has the same payoff structure as unit $1$ but with parameters
$$
 \tilde Q_{min}, \tilde Q_{max}, \tilde \alpha, \tilde  \beta, \tilde Q_d, \tilde \eta, \tilde c_{run}, \mbox{ and }  c_{low}.
$$
This problem can now be handled in the same way as above by introducing the payoff function
\begin{equation*}
f_3(Q,P,t) = \max_{\delta \in [0,1]} \left \{f_1 (\delta Q,P,t) + \tilde{f}_1 ((1-\delta) Q,P,t) \right \}.
\end{equation*}
and studying the resulting 4-state optimal switching problem with payoff functions $f_0, f_1, \tilde f_1$ and $f_3$, depending on whether no, either or both of the units are active. 
\end{remark}

\subsection{Switching costs}
Allowing for costs when switching between different states is one of the major benefits of an optimal switching approach to production planning. However, these costs (along with the risks involved when switching) are difficult to quantify, and may be so even for an experienced operator. In particular, they depend on the trait of optimization, e.g., maximizing profit, minimizing risk of stoppage, or minimizing wear of the units. Due to the uncertainty regarding the size of the switching costs, we study the outcome with different (constant) switching costs in Section \ref{sec:results} but refrain from explicitly modeling them in the current paper.

\subsection{Electricity spot price}
For computational simplicity, we have chosen to consider a deterministic electricity price $P$ in the current paper. To further facilitate the analysis and understanding of the results,	 we from here on in let $P_t\equiv P_0$ and perform our numerical analysis with constant electricity price. {However, we stress and reiterate that a time-varying \textit{deterministic} electricity price causes no other problems than obstructing the intuitive understanding of the results. }


\section{Solving the system of PDEs} \label{sec:numericsPDE}

With the parameters set above, we are ready to turn to solving \eqref{pde}. 
More precisely, we will solve a discretized version of \eqref{pde} based on the time-discretization 
$$\{t_0=0, t_1,t_2,\ldots,t_N=T\}$$
with $\Delta t = t_{n+1}-t_n= T/N$, and to do so we utilize an iterative implementation of the Crank-Nicolson scheme outlined in Algorithm \ref{alg:pseudo}. 

From here on in, we let $\{Q_{n}\equiv Q_{t_n}\}_{0\leq n \leq N}$  and $\{P_{n}\equiv P_{t_n}\}_{0\leq n \leq N}$ denote the discrete versions of $Q$ and $P$ based on this discretization. Moreover, starting at time $t_k$ with current flow $Q_{k}$, and given a discrete forecast $\{F_{t_{k+1}}, \dots, F_{t_{k+l}}\}$ of future flows we define the discrete versions of \eqref{eq:sigma} and \eqref{eq:b} by
$$
\sigma_k(Q_n, t_n)   = \sigma Q_n
$$
and
\begin{align} \label{eq:b2}
b_k(Q_{n}, t_n) &= 
\begin{cases}
b_k^{f} (t_n)							 &\textrm{if}\quad k< n \leq k+l +\ell,\\
 \left(r'_{t_n}+\frac{1}{2}\sigma^2-\kappa\left(\log Q_{n} - r_{t_n}\right)\right)\hat Q_{n} 	& \textrm{if}\quad k+l+\ell < n,
 \end{cases}
\end{align}
where 
$$
b_k^{f} (Q_{n},t_n)	 =  \left(g_k'({t_n})+\frac{1}{2}\sigma^2-\kappa\left(\log Q_{n} - g_k({t_n})\right)\right)Q_n,
$$
$g_k(t_n) := g_{t_k}(t_n)$ and $g'_{k}$ is calculated as a finite difference $g'_k(t_n) := \frac{g_k(t_n) - g_k(t_{n-1})}{\Delta t}$.


In the discretized version, and hence in the automated scheme to be constructed, the task is to find a strategy which maximizes the payoff
\begin{equation}\label{eq:discretisedpayoff}
\sum_{n=0}^{N-1} f_{\mu_{n}}(Q_{n},P_{n},{t_n}) \Delta t - \sum_{N-1 \geq i \geq 0} c_{\mu_{n-1}\mu_{n}}(Q_{n},P_{n},{t_n}),
\end{equation}
where $\mu_n \equiv \mu(t_n) : \{t_0,t_1,\ldots, t_N\} \to \{1,\ldots,m\}$ is a discrete valued function indicating the chosen production state between $t_n$ and $t_{n+1}$.

Recall that we are aiming at maximizing payoff from the current time $t_k$ up to some fixed terminal time $T$. Moreover, the function $b_k$ is updated as soon as our forecast is updated (in our case daily) and, consequently, from each starting time $t_k$ we will have a new function \eqref{eq:b2} and a new operator $\mathcal L^k$ in \eqref{pde} (based on the function $b_k$). Thus, our optimization problem differs slightly from time step to time step and we are bound to solve it repeatedly, for each point $t_k$ and with varying time span $[t_k,T]$. For $k$ fixed, our solution procedure is given by pseudo code of Algorithm \ref{alg:pseudo}.

\begin{algorithm} 

\KwResult{Numerical approximation of value function $u=\{u_1, \dots, u_m\}$}

Initialize $u$ by guessing a solution;

\While{totalerror $>$ totaltolerance}
{

		$u_{old} = u$;
		
	totalerror = 0;
	
	\For{time from $t_{N-1}$ to $t_k$ }{
		
		\While{spatialerror $>$ spatialtolerance}{
			
			spatialerror = 0;
			
			\For{all spatial points $q_j$}{
			
  				\For{$i \in \{1,\dots, m\}$}{
  				
 					Dummy variable $y$ = the Crank-Nicolson scheme for $u_i(q_j,t)$;
 					
  					$obstacle_i = \max_{j \neq i} \{u_j - c_{ij}\}$;
  					
  			    		$y = \max \{y,\text{$obstacle_i$} \}$;
  			    		
  			    		spatialerror = spatialerror + $|u_i(q_j,t)-y| $;
  			    		
  			    		$u_i(q_j,t)=y$;
				}
  			}
  			Interpolate linearly at the edges of the spatial grid
		}
	}
	totalerror = $|u_{old} - u|$;
 }
\caption{A numerical algorithm for estimating solutions to \eqref{pde}.}
\label{alg:pseudo}
\end{algorithm}

\section{Results}  \label{sec:results}

This section contains results of our parameter estimation and shows the performance of our PDE-based strategy. 

\subsection{Parameter values} \label{sec:parametervalues}
We test our model using flow data from the Swedish river Sävarån\footnote{The choice of Sävarån as our data source is made for no other reason than that the river lies close to Umeå University where most of the work on this paper was done.} during the years $1980-2018$, of which we use the first 35 years ($1980-2014$) for model calibration and the last four ($2015-2018$) for benchmarking. We also use data from $1980-2014$ for a long term evaluation of the performance of our scheme. Leap days are excluded in favour of a coherent presentation of the results. 

We study optimization over a one year horizon with possibility to change the state of production once every day, i.e.~$T=365$ days and $\Delta t= 1$ day. Our main reason for sticking with this coarse time-discretization is that the flow data available to us is given with one data point per day.

\begin{remark}
For hydropower plants with several units, i.e.,~plant II, the total flow of Sävarån is somewhat low and a non-regulated river with higher flow would be more suitable for our model. For simplicity of the presentation however, we have chosen to limit ourselves to a single river. Specific parameters, such as running costs, deterioration costs, etc., vary over time and between different hydropower plants and the parameters used here should, as with the choice of flow data, be seen merely as an example.
\end{remark}

\subsubsection*{Flow parameters}
\begin{figure}[h!]
\begin{subfigure}{.50\textwidth}
\centering
\includegraphics[width=1\linewidth]{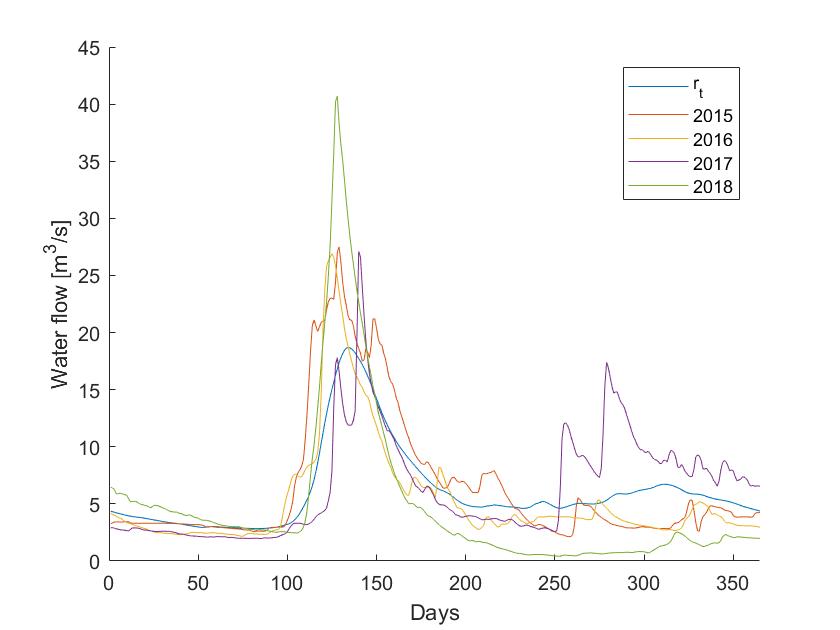}
\caption{The seasonal model $e^{r_t}$ together with actual flows $2015-2018$.}
\label{fig:flowmodel}
\end{subfigure}
\begin{subfigure}{.5\textwidth}
\centering
\includegraphics[width=1\linewidth]{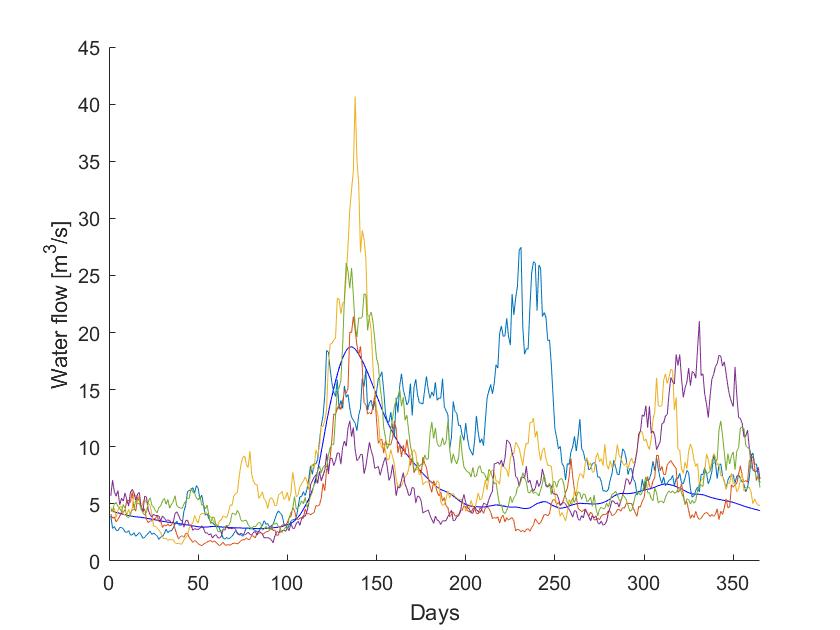}
\caption{The seasonal model $e^{r_t}$ together with 5 random realisations of \eqref{eq:sde}.}
\label{fig:simulatedflows}
\end{subfigure}
\caption{A visual comparison of the seasonal model $e^{r_t}$, solutions to \eqref{eq:sde} and actual flow during $2015-2018$.}
\end{figure}

From the historical flow data\footnote{Flow data were gathered from SMHI (http://vattenwebb.smhi.se/station) on June 15, 2020.} of Sävarån, the estimated values of the parameters in \eqref{eq:varacf} are $\kappa = 0.0208$ and $\sigma=0.1018$. Figure \ref{fig:flowmodel} shows $e^{r_t}$ and the flow during $2015-2018$ and Figure \ref{fig:simulatedflows} shows independent random realisations of \eqref{eq:sde} with these parameters.

\subsubsection*{Power plant parameters}

The parameters used in modeling our power plants are summarized and presented in Table \ref{tab:constants}.

\begin{figure}
\centering
\includegraphics[scale = 0.8]{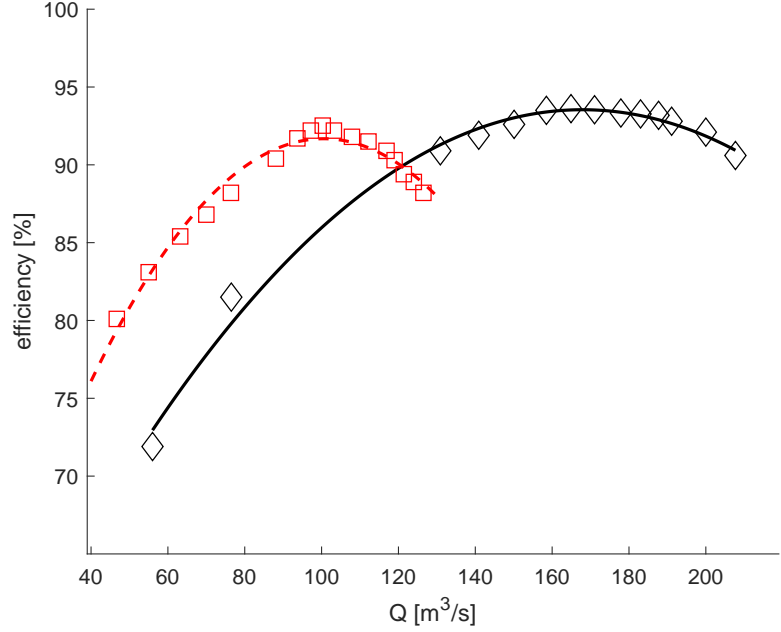}
\caption{Efficiency values of two Swedish Kaplan units with design flow of approximately 100 $m^3/s$ (red squares) and 170 $m^3/s$ (black diamonds).
The curves are $\eta(Q)$ defined in \eqref{eq:efficiency-curve-assumption} with least square fitted coefficients.
For the smaller unit $\alpha = 0.917$, $\beta = 0.430$ (red dotted) and for the larger unit  $\alpha = 0.935$, $\beta = 0.464$ (black solid).}
\label{fig:johej_eta}
\end{figure}

Starting with the efficiency curve of our power plant units, i.e., the coefficients of \eqref{eq:efficiency-curve-assumption}, Figure \ref{fig:johej_eta} shows measured efficiency of two Swedish Kaplan units together with the assumed efficiency curve in \eqref{eq:efficiency-curve-assumption} with least squared fitted coefficients $\alpha$ and $\beta$. Based on this data and the flow in Sävarån, we found it reasonable to choose the unit parameters as in Table \ref{tab:constants}.

The running cost is estimated from \cite{S12} to be approximately $1/5$ of the electricity price. The power of our unit is approximately $500$ $kW$, and with $P_t=P_0\equiv1$ it is thus reasonable to set $c_{run} = 100$~$m.u./h$. It is difficult to estimate the cost $c_{low}$ reflecting the running loss when the machine is run on too little water. The rotational speed of the turbine may drop so that the frequency of the produced electricity falls, resulting in a non-sellable production. We choose a value simply by multiplying $c_{run}$ with 10 and note that such a choice forces our algorithm to shut down the plant efficiently when the water supply fails.


It is not a trivial task to estimate a reasonable value of the switching costs. It may heavily depend on machine parameters related to the intake and specific properties of the turbine, tubes and the generator. Cost of personnel and environmental parameters such as local fish habitat may also be included, as well as the type of contract to which the electricity is sold. Due to these difficulties, we handle the switching costs as a parameter. In particular, we assume the switching costs to be constant and study the impact of varying this constant in Section \ref{sec:results}.

Moreover, when considering power plant 2, we assume that the cost of switching directly from state $0$ to state $2$ is cheaper than going through the intermediate state $1$, and vice versa. In particular, this implies that at any fixed time $t$, at most one switch is made. Exactly how much cheaper a direct switch from mode $0 \to 2$ (or vice versa) should be compared to the alternative $0 \to 1 \to 2$ depends on the actual power plant under consideration. Here, we simply assume that the alternative route is $50 \%$ more expensive. We summarize the switching costs in Table \ref{tab:swcosts}. Note that we have assumed all switching costs to be positive. When applicable, negative switching costs representing a gain rather than a cost when, e.g.,~reducing production capacity or moving to a more environmentally friendly production mode, can be used as well.

\begin{table}
\begin{center}
   \begin{tabular}{| l | l | l | l | }
    \hline
$P_0$         & 1 $[m.u. /kWh]$ &$c_{low}$      & $ 1000 $ $[m.u./h]$ \\  \hline
h           & 5  $[m]$        & $c_{run}$      & $ 100 $ $[m.u./h]$ \\ \hline
T           & 365 [days]      & $Q_{min}$      & 5 $[m^3/s]$  \\ \hline
$\alpha$    & 0.92            &$Q_{max}$      &13 $[m^3/s]$   \\ \hline
$\beta$     & 0.45            &   $Q_{d}$        &10 $[m^3/s]$    \\ \hline
\end{tabular}
\end{center}
\caption{Parameter values used in our numerical investigation.}
\label{tab:constants}
\end{table}

\begin{table}
\begin{center}
   \begin{tabular}{| l | c | c | c |}
    \hline
$c_{ij}  $ &   $j=   0$        &$j= 1$  &    $  j=   2 $  \\ \hline
$i=0  $     &$0$ &  $C $  & $1.5C$   \\  \hline
$i=1$       &$ C$ &  $0$& $C $\\  \hline
$i=2$       & $1.5C$ & $C$& $0$ \\ \hline
\end{tabular}
\end{center}
\caption{Relative switching costs. $C$ is a fixed constant determined in Section \ref{sec:results}.}
\label{tab:swcosts}
\end{table}
\subsubsection*{Forecasts}
\begin{figure}
\begin{subfigure}{.5\textwidth}
  \centering
  \includegraphics[width=1\linewidth]{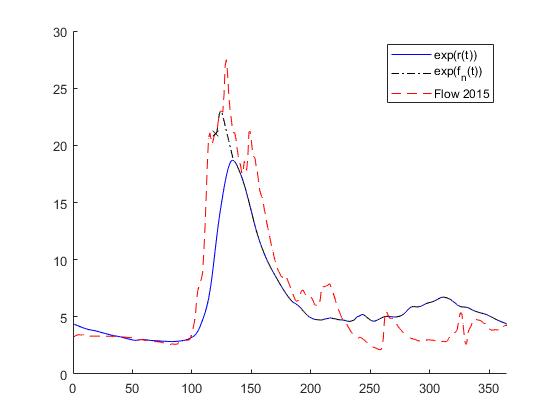}
  \caption{Plot of $e^{r(t)}$, the true flow for year 2015, \\and $e^{g_k(t)}$ for $k=120$, $l=5$ and $\ell=10$.}
  \label{fig:sub-first}
\end{subfigure} \hfill
\begin{subfigure}{.5\textwidth}
  \centering
  \includegraphics[width=1\linewidth]{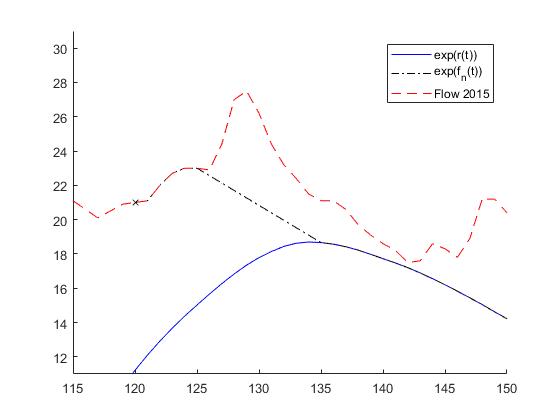}
  \caption{Zoom of figure (a)}
  \label{fig:sub-second}
\end{subfigure}

\begin{subfigure}{.5\textwidth}
  \centering
  \includegraphics[width=1\linewidth]{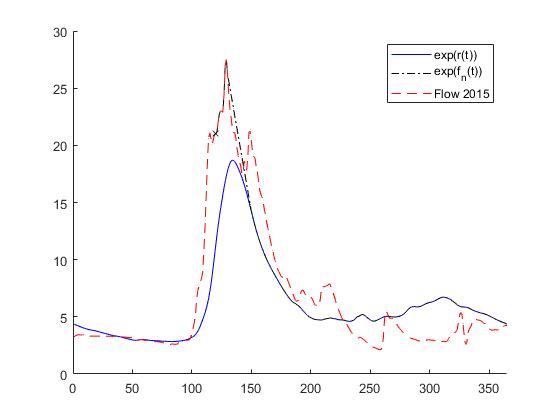}
  \caption{Plot of $e^{r(t)}$, the true flow for year 2015, \\and $e^{g_k(t)}$ for $k=120$, $l=10$ and $\ell=20$.}
  \label{fig:sub-first}
\end{subfigure}
\begin{subfigure}{.5\textwidth}
  \centering
  \includegraphics[width=1\linewidth]{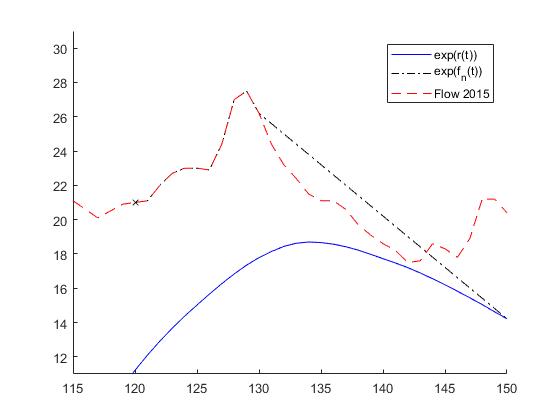}
  \caption{Zoom of figure (c)}
  \label{fig:sub-second}
\end{subfigure}
\caption{The impact of forecasts in relation to the long-term expected value $r_t$. Note that a stochastic term is also present in the flow model \eqref{eq:sde!} and hence the plotted figure above only corresponds to an estimate of the expected value of the future flow.}
\label{fig:forecasts}
\end{figure}

As already stated in the introduction, our main purpose is to highlight the use of optimal switching theory in production planning for ROR hydropower plants. The ambition is \textbf{not} to provide methods for forecasting river flow, and to avoid such discussions, we will simply use the true flow as forecast. However, we keep the stochastic component in the dynamics \eqref{eq:sde!} unchanged so that, even with forecast applied, our model does not know the future flow with certainty. Instead, this ``forecast'' only provides our model with a more accurate estimate of the average flow during the validity period of the forecast. We depict the impact of forecasts in Figure \ref{fig:forecasts}.

\subsection{Comparison of different strategies}
We will benchmark the performance of our strategy to the \textit{a fortiori} optimum, i.e.,~the optimal strategy in hindsight, with all information available. From a practical point of view, this value can only be achieved with certainty by ``looking into the future'', using the future flow of the river when making decisions. As this is of course impossible in practice, our comparison is slightly skew to the disadvantage of the model presented here. However, this value, the theoretical maximum output of the plant, is indeed achievable and should therefore be considered as the ultimate goal in any attempt to construct a production strategy. We emphasize and stress that, although the benchmark strategy can be found only in hindsight, our PDE-based strategy uses no other information when making a decision at time $t_k$ than historical information up to that point and the forecast starting at $t_k$. 

As a comparison, we also show the result of using a \textit{na\"{\i}ve strategy} in which the manager always switches to the production mode which momentarily has the highest payoff. To ease the presentation, we assume in all cases that the starting state is $0$, i.e.,~that the plant is ``off'' at the beginning of the planning period.

To get comparable results, we first calculate a fixed constant value $D$, depending on the power plant under consideration, but not on the switching costs or the flow of the river. More precisely, $D$ is calculated as the profit generated by the plant if it works at maximum capacity for a full year without interruptions (and starting in the most beneficial state), i.e.,
$$
D= f_i(Q_{max}, P_0, t) * 365,
$$ 
where $i =1$ or $i=2$ depending on the plant under consideration and $t \in [0,365]$ is arbitrary (since \eqref{eq:payoff} and \eqref{eq:f3} are independent of $t$). After determining $D$, the switching cost constant $C$ (cf. Table \ref{tab:swcosts}) is taken as a fixed percentage of $D$. 

Lastly, for comparison of strategies, we use the quotient $\gamma(\mu)$, defined as the total payoff from the strategy $\mu$ divided by $D$,
$$
\gamma(\mu) =   \frac{J_0(q_0,0,\mu)}{D},
$$
where $q_0$ is the flow at the starting time $t_0=0$ and $J_0$ is as defined in \eqref{eq:OSP}.

Results are provided for $T=365$ days with forecast lengths $l \in \{0,5,10\}$ days and with linear return to the long-term mean $e^{r_t}$ over $\ell=20$ days. Here, $l=0$ means no forecast. We give detailed descriptions of the suggested strategies for 2015 and show summarized results for $2016-2018$. Moreover, we provide results on the long term performance of our strategy by comparing results from the years $1980-2014$ for a fixed parameter set.
\footnote{It should be noted that, for convenience, this data set is the same as that used for calibrating parameters. Thus, for each year in the long term evaluation, the data tested is part of the data used for calibration. However, a single year out of the 35 used has minimal impact on the end calibration and the long term results are therefore still valid.}

\begin{figure}
\centering
\includegraphics[scale=0.6]{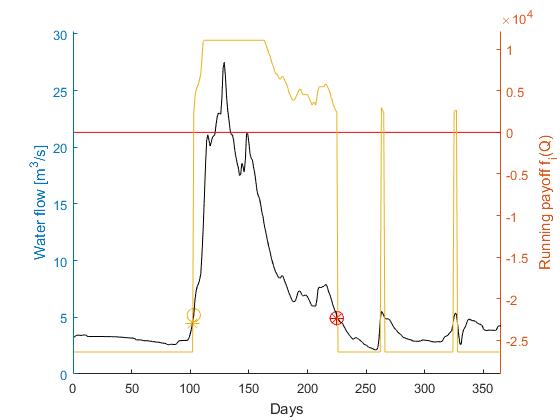}

\caption{Strategies for plant I during 2015 with $C/D=0.01$. The black line represents water flow and the yellow and red lines represent running payoff for states 1 and 0, respectively. Yellow circles (asterisks) represent the action of moving to state 1 and red circles (asterisks) the action of moving to state 0 under the PDE strategy with $l=10$ (optimal strategy).}
\label{fig:p1y1results}
\end{figure}

\begin{table}
 \centering
\begin{tabular}{l | l}
\textbf{Strategy}      & \textbf{(Day of action, Move to state)}   \\
\hline
Optimal                   &     (102,1) \quad (225,0)         \\    \hline
PDE, $l=0$               &     (103,1) \quad (226,0)          \\  \hline
PDE, $l=5$ &               (103,1) \quad (225,0)          \\  \hline
PDE, $l=10$              &  (103,1) \quad (225,0)    \\ \hline
Na\"ive             &    (102,1) \quad (225,0)    \quad  (262,1) \quad (265,0) \quad  (324,1) \quad (327,0)            \\ \hline
\end{tabular}
\caption{Strategies for plant I during 2015 with $C/D=0.01$.}
\label{tab:strategyp1}
\end{table}

\begin{figure}
\begin{subfigure}{.5\textwidth}
  \centering
  \includegraphics[width=1\linewidth]{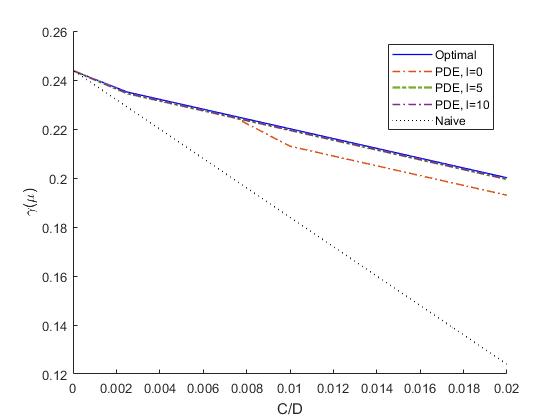}
  \caption{2015}
  \label{fig:sub-first}
\end{subfigure} \hfill
\begin{subfigure}{.5\textwidth}
  \centering
  \includegraphics[width=1\linewidth]{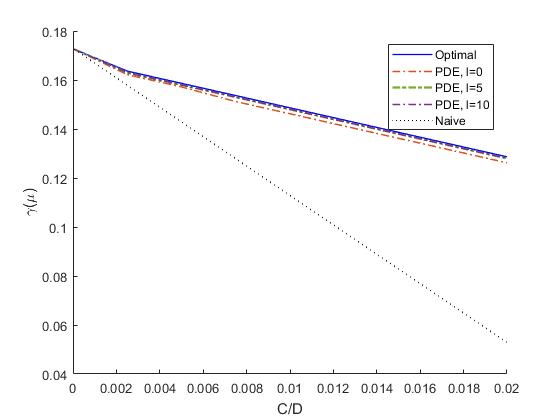}
  \caption{2016}
  \label{fig:sub-second}
\end{subfigure}

\begin{subfigure}{.5\textwidth}
  \centering
  \includegraphics[width=1\linewidth]{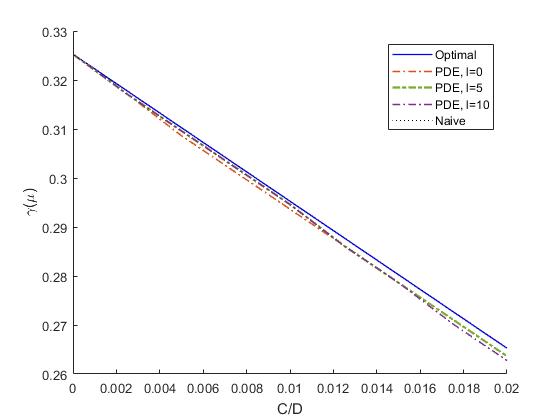}
  \caption{2017}
  \label{fig:sub-first}
\end{subfigure}
\begin{subfigure}{.5\textwidth}
  \centering
  \includegraphics[width=1\linewidth]{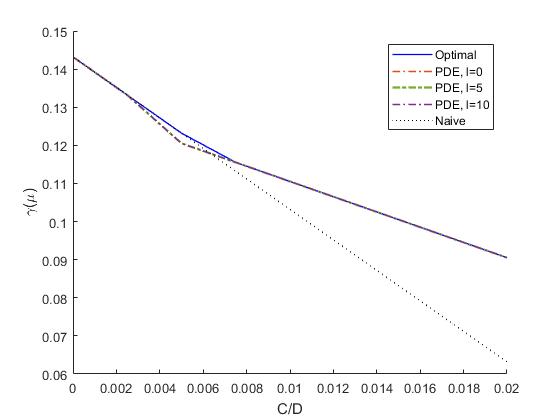}
 \caption{2018}
  \label{fig:sub-second}
\end{subfigure}
\caption{Relative payoff for plant I during $2015-2018$ as a function of $C/D$.}
\label{fig:relativesp1}
\end{figure}

\begin{figure}
\centering
  \includegraphics[scale=0.6]{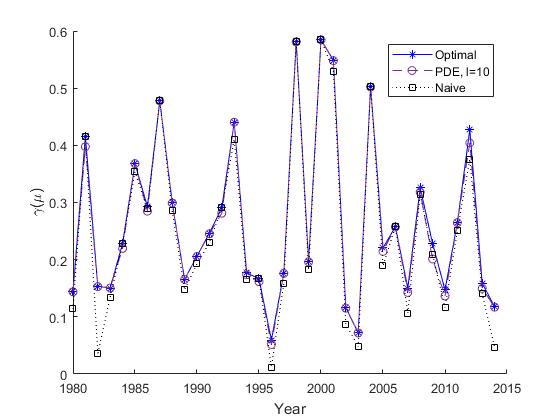}
  \caption{The PDE strategy consistently gives payoff close to the optimal.  Here, $C/D=0.01$ and $l=10$. The average quotient $\gamma$ for plant I over the years $1980-2014$ is $0.2674 $, $0.2628$, and $0.2466$ for the optimal strategy, the PDE strategy and the na\"{\i}ve strategy, respectively. (Lines are present for visual aid only.)} 
  \label{fig:averageresultsp1}
\end{figure}

\begin{figure}
\centering
\includegraphics[scale=0.6]{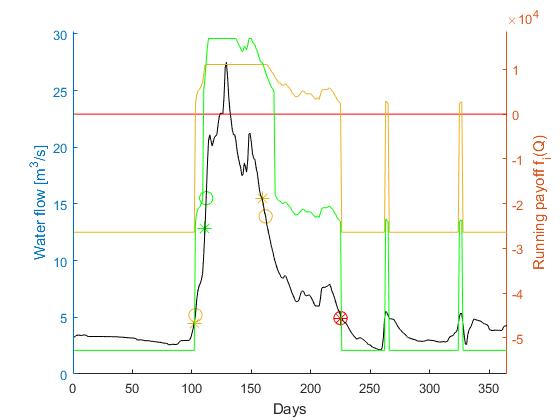}

\caption{Strategies for plant II during 2015 with $C/D=0.01$. The black line represents the water flow and the green, yellow and red lines represent running payoff for states 2, 1, and 0, respectively. Green circles (asterisks) represent the action of moving to state 2, yellow circles (asterisks) the action of moving to state 1, and red circles the action of moving to state 0 under the PDE strategy with $l=10$ (optimal strategy).}
\label{fig:p2y1results}
\end{figure}

\begin{table}
 \centering
\begin{tabular}{l | l}
\textbf{Strategy}      & \textbf{(Day of action, Move to state)}   \\
\hline
Optimal                   &     (102,1) \quad (111,2)   \quad (159,1) \quad     (225,0)       \\    \hline
PDE, $l=0$               &    (103,1) \quad (112,2)   \quad (162,1) \quad     (226,0)                        \\  \hline
PDE, $l=5$ &                   (103,1) \quad (112,2)   \quad (162,1) \quad     (225,0)              \\  \hline
PDE, $l=10$              &  (103,1) \quad (112,2)   \quad (162,1) \quad     (225,0)         \\ \hline
Na\"ive             &   (102,1) \quad (111,2)   \quad (159,1) \quad     (225,0) \quad (262,1) \quad (265,0)   \quad (324,1) \quad     (327,0)    \\ \hline
\end{tabular}
\caption{Strategies for plant II during 2015 with $C/D=0.01$.}
\label{tab:strategyp2}
\end{table}

\begin{figure}
\begin{subfigure}{.5\textwidth}
  \centering
  \includegraphics[width=1\linewidth]{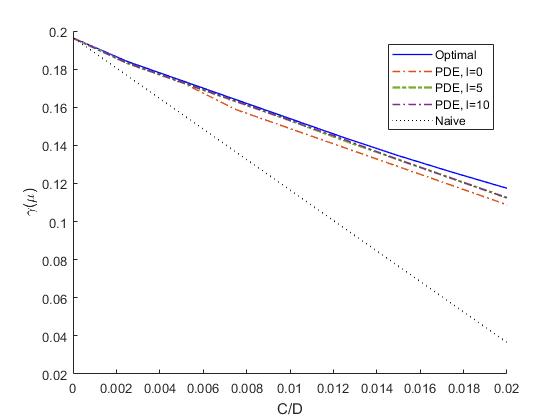}
  \caption{2015}
  \label{fig:sub-first}
\end{subfigure} \hfill
\begin{subfigure}{.5\textwidth}
  \centering
  \includegraphics[width=1\linewidth]{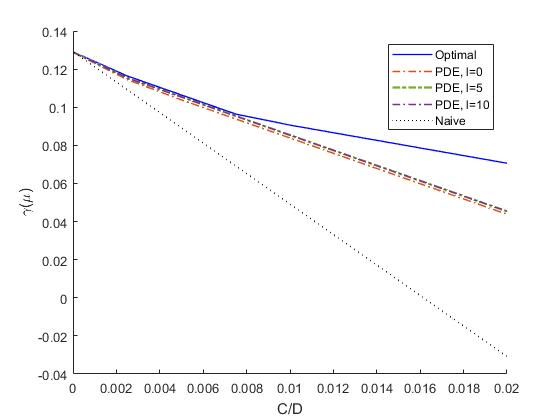}
  \caption{2016}
  \label{fig:sub-second}
\end{subfigure}

\begin{subfigure}{.5\textwidth}
  \centering
  \includegraphics[width=1\linewidth]{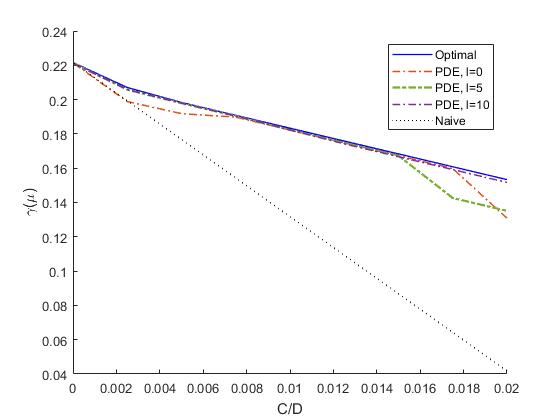}
  \caption{2017}
  \label{fig:sub-first}
\end{subfigure}
\begin{subfigure}{.5\textwidth}
  \centering
  \includegraphics[width=1\linewidth]{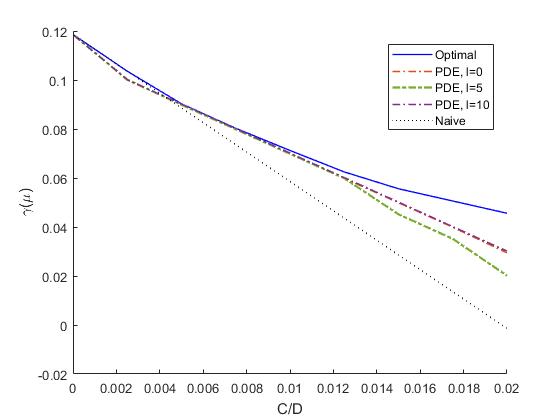}
  \caption{2018}
  \label{fig:sub-second}
\end{subfigure}
\caption{Relative payoff for plant II during $2015-2018$ as a function of $C/D$.}
\label{fig:relativesp2}
\end{figure}

\begin{figure}
\centering
  \includegraphics[scale=0.6]{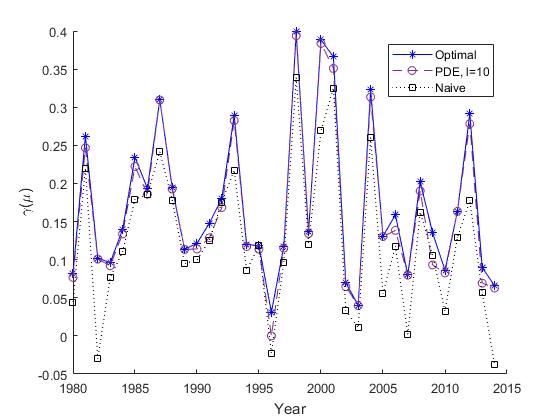}
  \caption{The PDE strategy consistently gives payoff close to the optimal. Here, $C/D=0.01$ and $l=10$. The average quotient $\gamma$ for plant II over the years $1980-2014$ is $0.1708 $, $0.1624$, and $0.1244$ for the optimal strategy, the PDE strategy and the na\"{\i}ve strategy, respectively. (Lines are present for visual aid only.)}
  \label{fig:averageresultsp2}
\end{figure}

The different production schemes (optimal, PDE-based, na\"{\i}ve) for plant I (plant II) during 2015, $l \in \{0,5,10\}$ and $C/D=0.01$ are presented in Figure \ref{fig:p1y1results} (Figure \ref{fig:p2y1results}) and Table \ref{tab:strategyp1} (Table \ref{tab:strategyp2}). The relative payoffs as a function of $C/D$ for $l \in \{0,5,10\}$ are given for all years in Figure \ref{fig:relativesp1} (Figure \ref{fig:relativesp2}). The long term performance of the PDE-based strategy for plant I (II) with $C/D=0.01$ is presented in Figure \ref{fig:averageresultsp1} (Figure \ref{fig:averageresultsp2}).

The PDE-based strategy in most cases performs very close to the optimal strategy, with a difference of less than 2 \% (5 \%) from the theoretical maximum in the long term tests for plant I (plant II) with $C/D=0.01$ and 10 days forecast. The most common mistake of the PDE-based strategy compared to the optimal is delaying the decision to change production mode. In all but a few cases, longer forecast also means better results. 


\section{Discussion} \label{sec:discussion}

The optimal switching theory is designed for maximizing the average payoff over a long period of time, but as our results show, it performs well also on single years. When comparing our strategy to the optimal one, we see that the differences between the strategies arise as our strategy occasionally recommends switching mode too late due to uncertainty regarding the future flow, see, e.g., Figure \ref{fig:p2y1results}. Most often, the decision is only late by one day and with a finer time discretization, these differences would most likely disappear or, at least, the discrepancy would be much smaller. We also see that by introducing forecasts, we are able to remove this gap entirely in many cases, see Figures \ref{fig:relativesp1} and \ref{fig:relativesp2}. 

Our (artificial) forecast includes uncertainty from the first forecasted day and reduction in this uncertainty, which is reasonable and possible by, e.g., upstream measurements, would also help removing any delay in the decision making. Indeed, in our SDE model, we assumed the uncertainty to be the same regardless of if a forecast was introduced or not. If the uncertainty of the forecast is known, one could introduce a new parameter $\sigma^f_k$ in a similar fashion as for $b^f_k$ and let the forecast influence also the stochastic volatility of the flow, having lower volatility close to the current time $t_k$ and increasing volatility further in the future. We have chosen not to alter the volatility $\sigma$ during the forecast period, partly to keep our model as simple as possible, and partly to avoid the need to construct forecasts (which is outside the scope of the current paper). 

Already without forecast, our PDE-based strategy outperforms the na\"{\i}ve strategy in most cases, even for small values of $C/D$, and in many cases also finds the truly optimal strategy or something very close to optimal. In the rare event that the na\"{\i}ve strategy performs as good or better than the PDE-based strategy, it is because the na\"{\i}ve strategy happens to be optimal. In these cases, the difference between the optimal and the PDE-based strategies is small. 

Our results show, as should be expected, that longer forecasts means better results. However, on a few occasions, this is not the case, see, e.g., Figure \ref{fig:relativesp1} (c), where, for $C/D \geq 0.016$, we perform worse with a 10 day forecast than with a 5 day forecast. The reason for this (rather unsatisfactory) result is that, when large fluctuations up and down in the flow happens in a just a few days time, a longer forecast may capture both directions of the movement whereas a shorter only captures one direction, triggering a decision to open/close if the current flow is at or close to a level at which different payoff functions intersect, whereas with longer forecast, the uncertainty introduced by forecasting also a rapid downward movement delays the decision slightly when the cost of switching is high. When run repeatedly over a large number of years, the decision made from the longer forecast would perform better \textit{on average}, but it may come up short in a single year. Luckily, as in the comparison with optimal na\"{\i}ve strategies, the deviation in the final payoff is small on these occasions.

Our model is calibrated to a constant electricity spot price $P_0$ for convenience when interpreting results and we repeat that a time-varying deterministic electricity price causes no problems other than parsing the results. However, our model could also be calibrated to a random price process $P_t$ as well. There is no (theoretical) restriction in the number of underlying Markovian It\^o processes our model can handle, so that allowing for such calibration is merely a question of computational power. In fact, not even the Markov property is a restriction as any discretized random process can be made Markovian at the cost of increasing its dimension. However, the cost of increasing the number of random sources is that the underlying optimization problem, which here is solved by PDE-methods, increases in dimensionality at the same rate. In practice, as long as the random sources are few, say 2 or 3, our approach based on numerical solutions of PDE can be used to find a solution. When the dimensionality increases even further, the PDE-methods become computationally heavy and other ways of attacking the resulting optimal switching problem may be preferable, e.g.~Monte Carlo-methods as in \cite{ACLP12,BGSS20}. In the current setting, the algorithm for obtaining our strategy is run in only a few minutes on a standard laptop computer and is thus more than sufficient for the purpose of the current paper.

\subsection{Concluding remarks and future research} \label{sec:concludingremarks}
In this paper we have, to the authors knowledge for the first time, used the mathematical optimal switching theory to create hydropower production plans which can incorporate random water flow and non-negligible costs of switching between different operational modes. Although our setup is somewhat simplified to keep the analysis of the results tractable, the results are satisfying, showing that automatic optimal switching schemes can perform close to the theoretical maximum already with small computational effort. Moreover, in our study the difference between our model and a na\"{\i}ve approach increased with the number of available production modes, indicating that the decision support provided by optimal switching theory becomes increasingly valuable as the complexity of the underlying problem increases.

An interesting theoretical continuation of the work initiated in this paper would be to study hydropower plants with a dam or hydropower plants of pumped-storage type. At the practitioners' level, a natural next step would be to adapt the current scheme to a real hydropower plant and to benchmark its performance to the production strategy currently in use. 

It is our firm believe that continued work in this direction, utilizing advanced stochastic control theory for production planning, will play an important role in making today's energy system more effective. Indeed, the theory outlined in this paper has potential to become as useful in the energy sector as the closely related theory of optimal stopping and singular- and non-singular control has become in financial mathematics. 

\section*{Acknowledgements}
The work of Niklas L. P. Lundstr\"om was partially supported by the Swedish research council grant 2018-03743.
Marcus Olofsson gratefully acknowledges the support from the Centre for Interdisciplinary Mathematics at Uppsala University.

\bibliographystyle{amsalpha}

\end{document}